\documentclass[a4paper, 12pt]{amsart}
\usepackage{mathrsfs,amssymb,amscd,amsmath}
\usepackage[T1]{fontenc}
\usepackage[utf8]{inputenc}
\usepackage[alphabetic]{amsrefs}
\usepackage[pdfauthor={Yiannis N. Petridis and Morten S. Risager},
             pdftitle={Higher order deformations of hyperbolic spectra},
             pdfsubject={Primary 11F72; Secondary  11M36, 34E10},
             pdfkeywords={},
             pdfproducer={LaTeX2e with hyperref},
             pdfcreator={Pdflatex},
             pdfdisplaydoctitle =true
 ]{hyperref}
\usepackage{pdfsync}
\usepackage{todonotes}

\newtheorem{theorem}{Theorem }[section]

\theoremstyle{definition}

\theoremstyle{remark}

\numberwithin{equation}{section}

\renewcommand{\Re}{\operatorname{Re}}
\renewcommand{\Im}{\operatorname{Im}}

\DeclareMathOperator*{\res}{res}

\def \G {{\Gamma}}
\def \g {{\gamma}}

\def \l {{\lambda}}

\def \R {{\mathbb R}}
\def \H {{\mathbb H}}

\def \C {{\mathbb C}}

\def \e {{\varepsilon }}
\def \GinfmodG {{\Gamma_{\!\!\infty}\!\!\setminus\!\Gamma}}
\def \GmodH {{\Gamma\setminus\H}}
\def \vol {\hbox{vol}}
\def \sl  {\hbox{SL}_2(\mathbb Z)}

\def \psl  {\hbox{PSL}_2(\mathbb R)}
\def \GmodH {\G\backslash\mathbb{H}}

\newcommand{\modsym}[2]{\langle #1 ,  #2  \rangle }

\newcommand{\norm}[1]{\left\lVert #1 \right\rVert}
\newcommand{\abs}[1]{\left\lvert #1 \right\rvert}
\newcommand{\inprod}[2]{\left \langle #1,#2 \right\rangle}
\newcommand{\tr}[1] {\hbox{tr}\left( #1\right)}

\dedicatory{In memory of Erik Balslev}
\title{Higher order deformations of hyperbolic spectra}
\author{Yiannis N. Petridis}
\address{Department of Mathematics, University College London, Gower Street, London WC1E 6BT, United Kingdom}
\email{i.petridis@ucl.ac.uk}

\author{Morten S. Risager}
\address{Department of Mathematical
  Sciences, University of Copenhagen, Universitetsparken 5, 2100
  Copenhagen \O , Denmark}
\email{risager@math.ku.dk}

\date{\today}

\begin{document}
\begin{abstract}
    This is an expanded writeup of a talk given by the second author at Erik Balslev's 75th birthday conference on October 1-2, 2010 at Aarhus University. We summarize our work on Fermi's golden rule and higher order phenomena for hyperbolic manifolds. A topic which occupied the last part of Erik Balslev's research.
\end{abstract}
\maketitle

\date{\today}
\section{Introduction}
In 1911 Herman Weyl \cite{Weyl:1911} proved  that the number $N(\lambda)$ of eigenvalues $\lambda_n$ less than $\lambda$ of the Dirichlet Laplacian on a bounded domain $X\subseteq \R^n$  with sufficiently nice boundary has the following asymptotic behaviour:
\begin{equation}
  N(\lambda)\sim \frac{\omega_d \vol{(X)}}{(2\pi)^d}\lambda ^{d/2}, \textrm{ as }\lambda\to\infty.
\end{equation}
Here $\omega_d$ is the volume of the unit ball in $\R^n$, and $\vol{(X)}$ is the volume of $X$. For us it is useful to know that Weyl's law holds for compact Riemannian manifolds, see \cite{MinakshisundaramPleijel:1949a, Buser:1992a}.
Weyl's law has been generalized and extended to many other cases, see e.g. \cite{Ivrii:2016c}.

In a seemingly unrelated direction Erich Hecke \cite{Hecke:1936} showed
that the zeta function $\zeta_\mathcal{O}$ of an \emph{imaginary} quadratic field $\mathcal{O}$ is related in a simple manner  to a certain modular form through Mellin transform. A modular form of weight $k$ is a holomorphic function on the upper half-plane $\H$ such that the differential $f(z)(dz)^{k/2}$ is invariant  under the action of certain subgroups of the full modular group $\sl$. Hans Maass \cite{Maass:1949a} investigated whether an analogous relation were true for \emph{real} quadratic fields. This led him to consider eigenfunctions of the hyperbolic Laplacian
\begin{equation*}
  -\Delta=-y^2\left(\frac{\partial^2}{\partial x^2}+\frac{\partial^2}{\partial y^2}\right)
\end{equation*}
on square integrable functions $f$ on the upper half-plane $\H$ transforming as $$f(\g z)=\chi(\g)f(z), \textrm{ for }\in \H, \gamma\in\Gamma.$$  Here $\G$ is a discrete subgroup of $\psl$ acting on $\H$ by linear fractional transformations, and $\chi:\G\to S^1$ is a unitary character. We denote the induced automorphic Laplacian by $L$. Such eigenfunctions have since been called Maass forms. Maass managed to show that -- at least for certain Hecke congruence groups $\G_0(N)$ and Dirichlet characters $\chi \mod N$ -- such forms exist. Moreover he showed that these forms are related to zeta functions of real quadratic fields in a way similar to how zeta functions of imaginary quadratic fields are related to modular forms. At that time it was not clear whether a single Maass form existed for $\sl$, $\chi=1$. On the other hand, Maass constructed non-holomorphic Eisenstein series, i.e. generalized eigenfunctions of $L$, but these are not square integrable.

Roelcke \cite{Roelcke:1953a, Roelcke:1966a} and Selberg \cite{Selberg:1956a, Selberg:1989a} gave a detailed description of the spectrum of $L$ when the hyperbolic volume $\vol{(\GmodH)}$ is finite. It consists of two parts:

\begin{enumerate}
\item A discrete set of eigenvalues
$$0\leq \lambda_0\leq \lambda_1 \leq \ldots\lambda_n\leq \ldots $$ This part may be finite or infinite, and does not have accumulation points.
\item   Furthermore, if $\GmodH$ is not compact, then the spectrum also contains a continuous part $[1/4, \infty[$ with multiplicity equal to the number of inequivalent open cusps for $(\Gamma, \chi)$.

\end{enumerate}
The continuous spectrum associated with the cusp $\mathfrak a$ is provided by Eisenstein series $E_{\mathfrak a}(z, s, \chi)$ for $s=1/2+it$, provided  the cusp $\mathfrak a$ is open, i.e. for its stabilizer $\Gamma_{\mathfrak a}$ in $\Gamma$ we have $\chi(\Gamma_{\mathfrak a})=1$.

We denote by $N_d(\l)=\#\{\l_n\leq \l\}$ the counting function for the \emph{discrete} part.
Using his newly developed trace formula, Selberg  proved the following groundbreaking result: if $\GmodH$ is not compact, but $\G$ is a congruence group and $\chi$ a Dirichlet character, then the set of discrete eigenvalues satisfy  Weyl's law, i.e.
\begin{equation}\label{weyls-law-discrete}
   N_d(\lambda)\sim\frac{\vol{(\GmodH)}}{4\pi}\lambda, \textrm{ as } \lambda\to\infty.
\end{equation}

Roelcke and Selberg independently speculated about the behaviour of  $N_d(\lambda)$ for general cofinite groups, i.e. $\G$ with $\vol{(\GmodH)}<\infty$, and general characters $\chi$. The belief that $N_d(\lambda)\to \infty$ as $\lambda \to\infty$ for such groups has been called the \emph{Roelcke--Selberg conjecture} by several authors, even if it is unclear
to what extent Roelcke and Selberg formally stated it as a conjecture. Even if there may be an infinite number of Maass forms,  they may or may not satisfy Weyl's law.

The difficulty of this conjecture lies in trying to count the discrete eigenvalues embedded in the continuous spectrum $[1/4, \infty[$. The general belief in the conjecture weakened after the work of Phillips and Sarnak on the stability of eigenvalues,  which we now describe, see also \cite{PhillipsSarnak:1985c, PhillipsSarnak:1985a, Sarnak:1990b, PhillipsSarnak:1992a}. A very  similar problem occurs in the study of the Schr\"odinger Hamiltonian for the helium atom, see \cite{Simon:1973a}.
In the physics literature embedded eigenvalues tend to be unstable and turn into scattering poles or resonances under perturbation. Resonances are poles of the analytic continuation of the resolvent in a second sheet (as opposed to the physical plane). The same phenomenon is true for the hyperbolic Laplacian. With the parametrization $\lambda=s(1-s)$ the second sheet corresponds to the left half-plane $\Re{s}<1/2$. The instability of embedded eigenvalues for the Schr\"odinger operator is described by Fermi's Golden Rule, proved rigorously by Simon in  \cite{Simon:1973a} following the work on analytic dilations of Balslev and Combes \cite{BalslevCombes:1971a}.  Phillips and Sarnak turned their attention to the analogous situation for hyperbolic surfaces with cusps.

Motivated by Selberg's trace formula, Phillips and Sarnak
defined the singular set. For a given eigenvalue $\lambda_j$  consider the  two values $s_j$ counting multiplicity satisfying $\lambda_j=s_j(1-s_j)$. \emph{The singular set} is then defined as follows: It is the multiset consisting of
\begin{enumerate}
  \item $s_j$ counted with the multiplicity of the corresponding eigenvalue.
  \item $\rho_j$ the poles of the scattering determinant $\varphi(s)$ counted with multiplicity the order of the pole.
  \item $1/2$ with multiplicity $(n+\tr{\Phi(1/2)})/2$ where $n$ is the number of open cusps of $(\G, \chi)$ and $\Phi$ is the scattering matrix related to $(\G,\chi)$.
\end{enumerate}
We refer to \cite{Selberg:1989a} for the definition of scattering matrix, open cusp, etc. Note that the above definition differs form \cite{PhillipsSarnak:1992a} by a rotation by $i$ followed by a shift of 1/2. Using the Lax--Phillips scattering theory \cite{LaxPhillips:1976a} as applied to automorphic functions, Phillips and Sarnak showed that the singular set is better behaved under deformations than the discrete spectrum.
We consider the following three types of deformations of $(\G, \chi)$:
\begin{enumerate}
\item  \label{character} Character deformations defined by
  \begin{align*}
      \chi_\varepsilon:&\G\to S^1\\
      &\g\mapsto \exp\left(2\pi i \varepsilon \int_{z_0}^{\gamma z_0}\alpha\right),
  \end{align*} where $\alpha=f(z)dz$ is a $\G$-invariant holomorphic 1-form, and $\varepsilon$ is a real parameter.
  \item \label{teichmuller} Real analytic deformations in Teichm\"uller space generated by $f$ a holomorphic cusp form of weight 4, see \cite{PhillipsSarnak:1985a} for details.
  \item \label{admissible} Real analytic compact deformations in the set of admissible surfaces, i.e.  Riemannian surfaces of finite area with hyperbolic ends, see \cite{Muller:1992a} for details.
\end{enumerate}

Phillips and Sarnak proved that in the cases \eqref{character},  \eqref{teichmuller}  the singular set has at most algebraic singularities. It follows in particular that, if $s(0)$ has multiplicity one, then $s(\varepsilon)$ is analytic for small $\varepsilon$. 
M\"uller \cite{Muller:1992a} extended this to case \eqref{admissible}, and Balslev \cite{Balslev:1997a} gave a different proof.

In all three cases described above the Laplacian $L(\varepsilon)$ admits a real analytic expansion
\begin{equation*}L(\varepsilon )=L(0)+\varepsilon L^{(1)}+\varepsilon^2\frac{L^{(2)}}{2}+\ldots
\end{equation*}
after possibly making a suitable conjugation, and adjustments of the corresponding metric high  in the cusps. See Section \ref{sec:character-def} for additional details.

Phillips and Sarnak identified a condition that will ensure that an embedded eigenvalue $\lambda_j=s_j(1-s_j)=1/4+t_j^2>1/4$ will dissolve into a resonance when $\varepsilon\neq 0$. For simplicity we restrict ourself to the case of only one open cusp. Let $E(z,1/2+it)$ be the generalized eigenfunction for the continuous spectrum at $1/4+t^2$ (see Section \ref{sec:Eisenstein}). Let furthermore $\hat s_j(\varepsilon)$ be the weighted mean of the branches of the singular points generated by splitting the eigenvalue $s_j(0)=s_j$ of multiplicity $m$ under perturbation, i.e.
\begin{equation*}
  \hat s_j(\varepsilon)=\frac{1}{m}\sum_{k=1}^ms_{j,k}(\varepsilon).
\end{equation*}
Let $u_{j,1}, \ldots u_{j,m}$ be an orthonormal basis of the eigenspace of $\lambda_j$.
\begin{theorem}[Fermi's Golden Rule]\label{FGR}
  If $\inprod{L^{(1)}u_{j,k}}{E(\cdot,1/2+it_j)}\neq 0$ for some $k$, then some of the  eigenvalues with eigenvalue $\lambda_j$   turn into  resonances under the perturbation. More precisely:
  \begin{equation*} \Re \hat s_j^{(2)}(0)=- \frac{1}{4t_j^2}\sum_{k=1}^m\abs{\inprod{L^{(1)}u_{j,k}}{E(\cdot,1/2+it_j)}}^2.
  \end{equation*}
\end{theorem}
For $m=1$ this is Eq. (5.29) in \cite{PhillipsSarnak:1992a}. For $m>1$ this is discussed in \cite{Petridis:1994c}.

Since the singular spectrum cannot move to the right under perturbation we always have $\Re \hat s_j^{(1)}(0)=0$, so $\Re \hat s_j^{(2)}(0)$ determines if $\hat s_j^{(2)}(0)$ moves to the left up to second order.

It turns out that for $\G$ a Hecke congruence group and perturbations of type \eqref{character} and \eqref{teichmuller} the dissolving condition \begin{equation}\label{Phillips-Sarnak}\inprod{L^{(1)}u_{j,k}}{E(\cdot,1/2+it_j)}\neq 0, \textrm{ for some }k=1,\ldots, l.\end{equation}
is equivalent to  the nonvanishing  of a special value of a Rankin--Selberg $L$-function. This allows one to use techniques from the analytic theory of $L$-functions to investigate how many eigenvalues are dissolved under perturbations. Luo \cite{Luo:2001a} succeeded in proving that a positive proportion of these special values of Rankin--Selberg $L$-functions are indeed non-zero. Assuming that multiplicities of the eigenvalues for the Laplacian of a fixed Hecke congruence groups are all bounded by the same common bound -- which is indeed expected to hold -- this allowed him to prove that a small deformation of the Hecke congruence groups does \emph{not} satisfy  Weyl's law, i.e. \eqref{weyls-law-discrete} does \emph{not} hold.

Phillips and Sarnak \cite{DeshouillersIwaniecPhillipsSarnak:1985a} conjectured something much stronger: the generic cofinite hyperbolic surface should only have finitely many discrete eigenvalues.

\subsection{Erik Balslev's interest in spectral deformations in the context of hyperbolic surfaces}
Throughout his career Erik Balslev was interested in various  properties of spectra of Schr\"o\-ding\-er operators. One of his major contributions to this field was the use of analytic dilation techniques in the setting of quantum mechanical many--body systems; see \cite{BalslevCombes:1971a}.

Balslev knew Ralph Phillips from his time in the United States in the 1960s and 1970s, and from Phillips numerous long-term visits to Denmark. Balslev was employed at Aarhus University, Denmark during most of his career.
In the autumn of 1991 Balslev was visiting Stanford and found himself in  Phillips'  office. Simultaneously Petridis went to the same office to  explain to Phillips his work on the genericity  of the $L^2$-eigenvalue $1/4$ and half-bound states, i.e. $E(z, 1/2)$, also known as nullvectors  \cite{Petridis:1994b}. This work was complementing that of Phillips and Sarnak.  Soon after this encounter  Balslev realized that the analytic dilation techniques, which he had used so effectively in the context of Schr\"odinger operators, could be used also in the context of deformations of hyperbolic surfaces. This realization led him to write \cite{Balslev:1997a}, where he reproved  much of the theory of Phillips and Sarnak, including Fermi's Golden rule, using analytic dilation techniques. 
It is well-known that the Eisenstein series $E(z, 1/2+it)$, which provide the continuous spectrum, have zero Fourier coefficient non-vanishing for $t\in\mathbb R^*$. However, the eigenfunctions with eigenvalue embedded in the continuous spectrum have vanishing zero Fourier coefficient.  Balslev introduced a family of operators $U(\lambda)$ acting only on the zero Fourier coefficient, corresponding to dilations in the hyperbolic distance for $\lambda $ real. For $\lambda $ complex the continuous spectrum $[1/4, \infty[$ of $\Delta$ is rotated by an angle $-2\arg \lambda$ to provide the continuous spectrum of the conjugated operator $U(\lambda)\Delta U(\lambda)^{-1}$. The embedded eigenvalues do not change location, so they become isolated. Resonances of $\Delta$ also turn into discrete eigenvalues (for appropriate choice of angle). This allowed him to use analytic perturbation theory and to reprove Fermi's Golden rule.

In the early 1990s Balslev met Alexei B. Venkov. This became the beginnning of a fruitful collaboration and close friendship which would last for the rest of Balslev's life. At this point Venkov had already been thinking about the Roelcke--Selberg conjecture and the Phillips--Sarnak conjecture for a long time \cite{Venkov:1979, Venkov:1990a, Venkov:1990c, Venkov:1992}, and together they started discussing the implications of Balslev's work \cite{Balslev:1997a}. Venkov started to visit Balslev at Aarhus University regularly, and in 2001 he joined their faculty. Together they worked on how to refine, use, and extend the deformation ideas of Phillips and Sarnak. This lead to several joint results on  Weyl's law \cite{BalslevVenkov:1998a, BalslevVenkov:2001a, BalslevVenkov:2007a} as well as on other related topics \cite{BalslevVenkov:2000, BalslevVenkov:2005a}.

\subsection{Higher order deformation}
The current work was inspired by the following question posed by  Erik Balslev to the authors: \emph{If the Phillips-Sarnak condition \eqref{Phillips-Sarnak} is \emph{not} satisfied, can one give simple conditions that ensures that an eigenvalue is dissolved?}
We will report on our work in this direction. We refer to \cite{PetridisRisager:2013a} for full details.
\bigskip

Understanding higher order deformations seems daunting at first. If one considers general expressions for the perturbation series of eigenvalues under analytic deformations one finds e.g. a 15-term expression for $\hat\lambda^{(4)}(0)$, see \cite[p. 80]{Kato:1976a}. We managed to find simpler expressions assuming that the lower order terms vanish, see Theorem \ref{maintheorem-for-characters} below.

Our motivation to understand what happens when the Phillips--Sarnak condition (\ref{Phillips-Sarnak}) is not satisfied came from the  numerical investigation by Farmer and Lemurell \cite{FarmerLemurell:2005a} and Avelin \cite{Avelin:2007b}. For a given cusp form Farmer and Lemurell found curves (branches) in Teichm\"uller space where a cusp form  for $\varepsilon=0$ remains cusp form, i.e. is not destroyed to \emph{any order}.  For specific even cusp forms Avelin  identified an analytic curve in Teichm\"uller space such that the movement of the poles of the scattering matrix gives a fourth order contact to the line $\Re(s)=1/2$. Our work aims to explain such phenomena theoretically.

\section{Stability of eigenvalues under character deformations}
We start by recalling a basic few properties of Eisenstein series. We refer to \cite{Iwaniec:2002a} for additional details.
\subsection{Standard non-hololorphic Eisenstein series}\label{sec:Eisenstein}
For simplicity of exposition we assume that $\Gamma$ has precisely one cusp, and that it is located at infinity. Placing the cusp at infinity can always be achieved by conjugation. Assume further that the stabiliser of the cusp is generated by $\gamma_\infty: z\mapsto z+1$, and that $\chi(\g_\infty)=1$. Recall that the standard non-holomorphic Eisenstein series is defined by
\begin{equation*}
  E(z, s,\chi)=\sum_{\gamma\in\GinfmodG}\overline\chi(\gamma)\Im(\gamma z)^s, \textrm{ when }\Re(s)>1,
\end{equation*}
and that it admits meromorphic continuation to $s\in \C$. Recall also that $E(z,1/2+it,\chi)$ is a generalized eigenfunction of $\Delta$ with eigenvalue $1/4+t^2$. These numbers, with $t\geq 0$ span the continuous spectrum. Clearly $E(\g z, s,\chi)=\chi(\g)E(z, s,\chi)$. The zero Fourier coefficient of $E(z, s,\chi)$ has the form
 \begin{equation*}
   \int_0^{1}E(z, s,\chi)dx=y^s+\varphi(s , \chi)y^{1-s}.
 \end{equation*}
This defines the scattering matrix $\varphi(s,\chi)$, which in the one cusp case is just a function. We recall that $\varphi(s,\chi)$ satisfies the functional equation
\begin{equation}
  \varphi(s, \chi)\varphi(1-s, \chi)=1,
\end{equation} and that it is unitary on the line $\Re(s)=1/2$.
Furthermore \begin{equation}\label{functional}
E(z, s,\chi)=\varphi(s, \chi)E(z, 1-s,\chi).
\end{equation}
We define
\begin{equation}M(T)=-\frac{1}{4\pi}\int_{-T}^T\frac{\varphi'}{\varphi}(1/2+it, \chi)dt.
\end{equation}
 Selberg proved that for \emph{all} cofinite groups $\G$ we have
\begin{equation*}
N_d(\lambda)+M(\sqrt{\lambda-1/4})\sim \frac{\vol{(\GmodH)}}{4\pi}\lambda.
\end{equation*}
Hence $N_d$ satisfies  Weyl's law precisely if
\begin{equation}\label{bound-M}
  M(T)=o(T^2).
\end{equation} Selberg proved that \eqref{bound-M} holds for $\G$ a congruence group. For such a group the scattering determinant $\varphi$ can be computed explicitly in terms of completed $L$-functions, and the bound \eqref{bound-M} follows from classical bounds on these $L$-functions.

\subsection{Character deformations}\label{sec:character-def}
For the rest of the paper we consider for simplicity  the case of character deformations. The cases of real analytic Teichm\"uller deformations and real analytic compact deformations within the set of admissible surfaces can be dealt with in a similar way. We refer to \cite{PetridisRisager:2013a} for additional details.

 Let $f$ be a cusp form of weight two, i.e. $f:\H\to\C$ is a holomorphic function that satisfies
 \begin{equation*}
   f(\gamma z)=(cz+d)^2f(z) \textrm { for all } \g\in \G,
\end{equation*}
and admits a Fourier expansion
\begin{equation}\label{Fourier-expansion}
  f(z)=\sum_{n=1}^\infty a_n e^{2\pi i n z}.
\end{equation}
Let $\omega=\Re(f(z)dz)$ be the corresponding real invariant 1-form, and let $\alpha$ be a compactly supported $1$-form in the same cohomology class as $\omega$.  We now define the modular symbol to be
\begin{equation*}
  \modsym{\gamma}{\alpha}:=-2\pi i\int_{z_0}^{\gamma{z_0}}\alpha.
\end{equation*}
Here $z_0$ is any point in $\H\cup\{i \infty \}$. The modular symbol is independent of the choice of $z_0$, the choice of compact form $\alpha$ in the same cohomology class,  as well as choice of path between $z_0$ and $\gamma z_0$. The modular symbol mapping  is an additive homomorphism, i.e. $\modsym{\gamma_1\gamma_2}{\alpha}=\modsym{\gamma_1}{\alpha}+\modsym{\gamma_2}{\alpha}$, and, moreover, vanishes at parabolic elements: $\modsym{\g_\infty}{\alpha}=0$.

With the help of modular symbols we create a one-parameter family of unitary characters. Consider now the unitary characters $$\chi_\varepsilon(\gamma)=\exp(\varepsilon\modsym{\g}{\alpha})$$ and the space
\[ L^{2}(\GmodH, \overline\chi_\varepsilon))=\left\{f:\H\to\C:f(\gamma z)=\overline\chi_\varepsilon(\g)f(z), \int_{\GmodH}\abs{f(z)}^2d\mu(z)<\infty
\right\}.\] Here $d\mu(z)=y^{-2}dxdy$ is the $\psl$-invariant measure on $\H$. We denote the induced automorphic Laplacian by $\tilde L(\varepsilon)$. We now conjugate this family of operators to the fixed space $L^2(\GmodH)$ by using  unitary operators
\begin{align*}
  U(\varepsilon): L^2(\GmodH) \to& L^{2}(\GmodH, \overline\chi_\varepsilon))\\
   f(z)\mapsto & \exp\left(2\pi i \varepsilon \int_{z_0}^z \alpha\right)f(z),
\end{align*}
and let $L(\varepsilon)=U^{-1}(\varepsilon)\tilde L(\varepsilon) U(\varepsilon)$. This new family of operators has the advantage of being defined on a fixed space. It is now a straightforward computation to show that on smooth functions $h$
\begin{equation*}
L(\varepsilon)h=\Delta h +\varepsilon L^{(1)}h+\frac{\varepsilon^2}{2}L^{(2)}h,
\end{equation*}
where
\begin{align}
  \label{infinitesimal}
\begin{split}
  L^{(1)}h&=4\pi i \inprod{dh}{\alpha}-2\pi i \delta(\alpha)h,\\
  L^{(2)}h&=-8\pi^2\inprod{\alpha}{\alpha}h.
\end{split}
\end{align}
Here
\begin{align*}
  \inprod{f_1dz+f_2d\overline z}{g_1dz+g_2d\overline z}&=2y^2\left(f_1\overline g_1+f_2\overline g_2\right),\\
\delta(pdx+qdy)&=-y^2(p_z+q_y).
\end{align*}
 We want to investigate whether embedded eigenvalues are destroyed under this perturbation, i.e. they turn into resonances.

\subsection{Goldfeld Eisenstein series} Goldfeld \cite{Goldfeld:1999a, Goldfeld:1999b} introduced  in the late 1990s  a generalization of the standard Eisenstein series $E(z,s)$, which has since been studied by several people.   It turns out that the stability of eigenvalues under perturbations as described above can be analyzed using such series.

The Goldfeld Eisenstein series, also known as Eisenstein series twisted by modular symbols, is defined by
\begin{equation*}
  E^n(z,s)=\sum_{\gamma\in\GinfmodG}\modsym{\g}{\alpha}^{n}\Im(\gamma z)^s, \textrm{ when }\Re(s)>1.
\end{equation*}
It is well establised \cite{PetridisRisager:2004a, JorgensonOSullivan:2008a} that $E^n(z,s)$ admits meromorphic continuation to $s\in \C$ . 

For $n>0$ the function $E^n(z,s)$ is not invariant but satisfies an $n$th order automorphy relation, i.e. $E^n(z,s)\in A^n_\G$ where $A^n_\G$ is defined recursively as follows: the set $A^0_\Gamma$ is simply the set of $\Gamma$-invariant functions on $\H$, and $A^n_\G$ consists on functions $f$ on $\H$ satisfying  $f(\g z)-f(z)\in A_\Gamma^{n-1}$ for all $\gamma\in \G$. For details on higher order Maass forms see e.g. \cite{BruggemanDiamantis:2012a}.

As $\Im (\gamma z)^s$ is formally an eigenfunction of $\Delta$, we have furthermore
\begin{equation*}(\Delta +s(1-s))E^{n}(z,s)=0.
\end{equation*}

There is a related \emph{invariant} function constructed by automorphizing $\left(-2\pi i \int_{z_0}^{ z}\alpha \right)^{n}\Im( z)^s$ as follows:
\begin{equation*}
  D^n(z,s)=\sum_{\gamma\in\GinfmodG}\left(-2\pi i \int_{z_0}^{\gamma z}\alpha \right)^{n}\Im(\gamma z)^s, \textrm{ when }\Re(s)>1.
\end{equation*}
Similarly to $E^n(z,s)$ the function $D^n(z,s)$ also admits meromorphic continuation to $s\in \C$. Indeed there is a simple way to relate the two:
\begin{equation*}
  D^{n}(z,s)=\sum_{j=0}^n\binom{n}{j}\left(-2\pi i \int_{z_0}^z\alpha \right)^{n-j}E^{j}(z,s).
\end{equation*}
The function $D^{n}(z,s)$ is not an eigenfunction of $\Delta$ but satisfies
\begin{equation}\label{laplacian-to-D}(\Delta +s(1-s))D^{n}(z,s)=-\binom{n}{1}L^{(1)}D^{n-1}(z,s)-\binom{n}{2}L^{(2)}D^{n-2}(z,s),
\end{equation}
where $L^{(1)}$ and $L^{(2)}$ are as in \eqref{infinitesimal}. Here we interpret $D^{n}(z,s)=0$, if $n$ is negative.

\subsection{Higher order Fermi's Golden Rules}

We are now ready to formulate our main theorem, which answers the question of Balslev mentioned in the introduction:
\begin{theorem}[\cite{PetridisRisager:2013a}]
  \label{maintheorem-for-characters}
  Let $s_j$ be a cuspidal eigenvalue for $L(0)$ and let $\hat s_j(\varepsilon)$ be the weighted mean of the branches of the singular points for $L(\varepsilon)$ generated by $s_j$. Assume that $$\hat s_j^{(l)}(0)=0, \quad \textrm{ for }l\leq 2(n-1).$$
  Then
\begin{enumerate}
\item $\hat s_j^{(2n-1)}(0)=0$,
\item
  $D^{n}(z, s)$ has at most a first order pole at $s_j$, and
 \item
 $\displaystyle
 \Re {\hat{s}}_j^{(2n)}(0)=-\frac{1}{2n}\binom{2n}{n}\norm{     \res_{s=s_j}D^n(z, s)}^2.$
\end{enumerate}
\end{theorem}

We  note that if
\begin{equation}\res_{s=s_j}D^{n}(z, s)\neq 0,
\end{equation}
 then at least one eigenvalue $s_j$ will be become a resonance under the deformation, so we may interpret this as a \emph{higher order vanishing condition}.

We also note that when $n=1$ this reduces to Theorem \ref{FGR}.  We note  that the Phillips--Sarnak vanishing condition can be formulated as
\begin{equation}\label{PS-reformulated}
\res_{s=s_j}D^{1}(z,s)\neq 0,
\end{equation}
since we have
\begin{equation*}\res_{s=s_j}D^{1}(z,s) =\sum_{k=1}^m \frac{c}{\pi^{s_j}}L(u_{j,k}\otimes f, s_j+1/2)\Gamma (s_j-1/2)u_{j, k}(z)   ,
\end{equation*} see \cite[Eq. (1.13)]{Petridis:2002a}  and \cite{PhillipsSarnak:1985a}.
\section{Relation to special values of Dirichlet series}
We first need to setup some additional notation.
An eigenfunction $u_j$ with eigenvalue $s_j(1-s_j)>1/4$ has a Fourier expansion
\begin{equation*}
  u_j(z)=\sum_{n\neq 0}b_n\sqrt{y}K_{s_j-1/2}(2\pi\abs{n}y)e^{2\pi i n x}
\end{equation*}
where $K_s(t)$ is the McDonald--Bessel function. We assume that $u_j$ has been normalized to have $L^2$-norm equal to 1.

The Phillips--Sarnak condition for dissolving a cuspidal eigenvalue \eqref{Phillips-Sarnak} can be expressed as the non-vanishing of a special value of a Rankin--Selberg $L$-function, see \cite{Sarnak:1990b, PhillipsSarnak:1985a}. This is defined as
\begin{equation*}
  L(u_j\otimes f,s)=\sum_{n=1}^\infty\frac{a_{n}b_{-n}}{n^{s+1/2}}\textrm{ for }\Re(s)>1.
\end{equation*}
 We will now explain that something similar happens for the higher order dissolving conditions in Theorem \ref{maintheorem-for-characters}.

We consider  the antiderivative of the cusp form $f$ inducing $\chi_\varepsilon$, i.e.
\begin{equation*}
  F(z)=\int_{i\infty}^{z} f(w)dw=\sum_{n=1}^\infty\frac{a_n}{2\pi i n}e^{2\pi i n z} .
\end{equation*}
Then for $\Re(s)>1$ we define the convergent Dirichlet series
\begin{equation}
  L(u_j\otimes F^2,s)=\sum_{n=1}^\infty\sum_{k_1+k_2=n}\frac{a_{k_1}}{k_1}\frac{a_{k_2}}{k_2}\frac{b_{-n}}{n^{s-1/2}}
\end{equation}
One can show that $L(u_j\otimes F^2,s)$ admits meromorphic continuation to $s\in \C$ and satisfies a functional equation relating its value at $s$ and $1-s$. The possible poles of $L(u_j\otimes F^2,s)$  are at  the singular points.

A holomorphic form of weight 2 as in Section \ref{sec:character-def} gives rise to two character deformations, namely those induced from $\omega_1=\Re(f(z)dz)$ and $\omega_2=\Re(if(z)dz)$.
\begin{theorem} Assume that the Phillips--Sarnak condition \eqref{PS-reformulated} at a cuspidal eigenvalue $s_j$ is \emph{not} satisfied for either of $\omega_1$, $\omega_2$. Then
$L(u_j\otimes F^2,s)$ has a removable singularity at $s_j$.

Assume further that $L(u_j\otimes F^2,s_j)\neq 0$. Then in all directions $\omega$ in the real span of $\omega_1, \omega_2$ with at most two exceptions we have
\begin{equation*}
\Re\hat s_j^{(4)}(0)\neq 0.
\end{equation*}
In particular there exists a cusp form with eigenvalue $s_j(1-s_j)$ that is dissolved in this direction.
\end{theorem}

We refer to \cite[Sec. 4.3]{PetridisRisager:2013a} for proofs of this theorem.
The function $L(u_j\otimes F^2,s)$ is not as well studied as the Rankin--Selberg $L$-function, and, although it does share many of its properties (continuation to $s\in \C$, functional equation, bounds on vertical lines), there are important differences. Most importantly $L(u_j\otimes F^2,s)$ does not admit an Euler product.

\section{Idea of proof}
We now indicate the main steps of Theorem~\ref{maintheorem-for-characters}, and refer to \cite{PetridisRisager:2013a} for details.

For the fixed group $\Gamma$ and the family of characters $\chi_{\varepsilon}$ we consider the scattering matrix $\varphi(s, \varepsilon)$. Besides properties that we have already stated one can show  that
\begin{equation}\label{symmetry}
  \varphi(s,\varepsilon)=\overline \varphi(\overline s, \varepsilon),
\end{equation}
see \cite[page 218, Remark 61]{Hejhal:1983a}.

We track  the movement of singular set close to an embedded eigenvalue $s_j(1-s_j)>1/4$ i.e. the embedded eigenvalue/resonance in the half-plane left of $\Re(s)=1/2$ using complex analysis, in particular, a simple variation of the argument principle. Define $\Lambda$ as the half circle $\gamma_1 (t)=ue^{it}+s_j$, $\pi/2\le t\le 3\pi/2$ followed by the line $\gamma_2 (t)=s_j+it$, $ -u\le t\le u$. Here   $u$ is chosen small enough, so that the only singular point for $\e=0$ inside  the ball $B(s_j, u)$ is $s_j$ with multiplicity $m=m(s_j)$. This contour is traversed counterclockwise.
For $\e$ sufficiently small the total multiplicities of  the singular points $s_j(\e)$ inside $B(s_j, u)$ is $m(s_j)$.

 We have
\begin{equation}\label{blah}m(\hat{s} (\e)-s_j)=-\frac{1}{2\pi i}\int_{\Lambda} (s-s_j)\frac{\varphi' (s, \e)}{\varphi (s, \e)}\, ds+\sum_{j\in C}(s_j(\e)-s_j),\end{equation}
where $C$ is indexing the cusp forms eigenbranches inside $B(s_j, u)$, i.e. the cusp forms that remain cusp forms. Let the last sum be denoted by $p(\e)$. The reason for using $\Lambda$ and not the whole $\partial B(s_j, u)$ is that on the right half-disc $\varphi (s,\e)$ has zeros, which we do not want to count. Note that by well-known properties of $\varphi(s,\e)$  \cite[Chapter 6]{Iwaniec:2002a} it has no zeroes in $\Lambda$.
Notice that $\overline{\int_{\gamma}f(s)\, ds}=\int_{\bar \gamma} \bar
f(\bar s)\, ds$ and, therefore we find  by  using (\ref{symmetry}) that
\begin{align*}m(\overline{\hat{s} (\e)-s_j})&=\frac{1}{2\pi i}\int_{\bar \Lambda} (s-\bar s_j)\overline{\left(\frac{\varphi' (\bar s, \e)}{\varphi (\bar s, \e)}\right)}\, ds+\overline{p(\e)}\\
  &=\frac{1}{2\pi i}\int_{\bar \Lambda}(s-\bar s_j)\frac{\varphi' (s, \e)}{\varphi (s, \e)}\, ds+\overline{p(\e)}.
\end{align*}
Denoting by $-\gamma$ the contour $\gamma$ traversed in the opposite direction, we get
\begin{align*}m(\overline{\hat{s} (\e) -s_j})&=-\frac{1}{2\pi i }\int_{-\bar \Lambda}(s-\bar s_j)\frac{\varphi'(s, \e)}{\varphi (s, \e)}\, ds+\overline{p(\e)}\\&=-\frac{1}{2\pi i}\int_{T^{-1}(-\bar\Lambda)}(1-w-\bar s_j)\frac{\varphi'(1-w, \e)}{\varphi (1-w, \e)}\, (-dw)+\overline{p(\e)},\end{align*}
where $s=T(w)=1-w$ is a conformal map. By \eqref{functional} we get
\[\varphi'(s, \e)\varphi (s, \e)-\varphi (s, \e)\varphi'(1-s, \e)=0,\] giving that \[\frac{\varphi'(s, \e)}{\varphi (s, \e)}=\frac{\varphi'(1-s, \e)}{\varphi(1-s, \e)}.\]
We plug this into the expression for $m(\overline{\hat{s} (\e) -s_j})$ to get
\begin{equation}\label{blahblahblah}m(\overline{\hat{s} (\e) -s_j})=-\frac{1}{2\pi i}\int_{T^{-1}(-\bar\Lambda)}(w- s_j)\frac{\varphi'(w, \e)}{\varphi (w, \e)}\, dw+\overline{p(\e)}. \end{equation}
We sum (\ref{blah}) and (\ref{blahblahblah}) and notice that the
cuspidal branch contributions cancel, because for a cuspidal branch
$s_{j,l}(\e)$ the function  $s_{j,l}(\e)-s_j$ is purely imaginary.
 We therefore conclude that
 \begin{eqnarray}\label{flocke} \nonumber 2m\Re (\hat{s}(\e)-s_j)&=&-\frac{1}{2\pi i}\int_{\Lambda+T^{-1}(-\bar \Lambda)}(s-s_j)\frac{\varphi'(s, \e)}{\varphi(s, \e)}\, ds\\&=&  -\frac{1}{2\pi i}\int_{\partial B(s_j,
 u)} (s-s_j)\frac{\varphi'(s, \e)}{\varphi(s, \e)}\, ds,\end{eqnarray}
 since the contribution from the line
 segment on $\Re(s)=1/2$ from $\Lambda$ and $T^{-1}(-\bar \Lambda)$ cancel.
  By uniform convergence
 we can differentiate the last formula in $\e$. We get
 \begin{align}\label{mosedebatteri}&\left.2m\frac{d^{2n}}{d\e^{2n}}\Re (\hat{s} (\e))\right|_{\e=0}=-\frac{1}{2\pi i}\int_{\partial B(s_j, u)}(s-s_j)\frac{d^{2n}}{d\e^{2n}}\left.\left(\frac{\varphi'(s, \e)}{\varphi (s, \e)}\right)\right|_{\e=0}\, ds\\
\nonumber &=-\frac{1}{2\pi i}\int_{\partial B(s_j, u)}(s-s_j)\sum_{k=0}^{2n}\binom{2n}{k} \left.\frac{d^{k}\varphi'(s, \e)}{d\e^{k}}\right|_{\e=0}\left.\frac{d^{2n-k}(\varphi (s, \e)^{-1})}{d\e^{2n-k}}\right|_{\e=0}\, ds.
\end{align}
 We can interchange the order of differentiation: 
 \[\frac{d^k}{d\e^k}\varphi'(s, \e)=\frac{d}{ds}\varphi^{(k)}(s).\]
 Note that the prime denotes derivative in $s$, whereas the $\varphi^{(k)} (s)$ denotes the $k$th derivative in $\varepsilon$ evaluated at $\varepsilon=0$. By differentiating $m$ times $\varphi(s,\e)^{-1}\varphi(s,\e)=1$ we find
 \begin{equation*}
   \sum_{k=0}^m\binom{m}{k}\left.\frac{d^{k}}{d\e^{k}}\varphi(s,\e)^{-1}\right|_{\e=0}\varphi^{(m-k)}(s,0)=0.
 \end{equation*}

These observations combined with \eqref{mosedebatteri} show that in order to compute $\left.2m\frac{d^{2n}}{d\e^{2n}}\Re (\hat{s} (\e))\right|_{\e=0}$ it suffices to understand the analytic behaviour of $\varphi^{(k)}(s)$ at $s=s_j$.

The general functional equation \eqref{functional} implies that
\begin{equation*}D^{n}(z,s)=\sum_{k=0}^n\binom{k}{n}\varphi^{(k)}(s)D^{n-k}(z,1-s).\end{equation*}
Combining this with \eqref{laplacian-to-D} and properties of the resolvent kernel $R(s)=(\Delta+s(1-s))^{-1}$ we can show that
\begin{equation}\label{phi-n-expr}\varphi^{(n)}(s)=\frac{1}{2s-1}\int_{\GmodH}E(z,s)\left(\binom{n}{1}L^{(1)}D^{n-1}(z,s)+\binom{n}{2}L^{(2)}D^{n-2}(z,s)\right)d\mu(z).\end{equation}
After some computations this and \eqref{mosedebatteri} lead to

\begin{align*}
  \left.2m\frac{d^{2n}}{d\e^{2n}}\Re (\hat{s}
    (\e))\right|_{\e=0}&=\frac{\res_{s=s_j}\varphi^{(2n)}(s,0)}{\varphi(s_j,0)}.
\end{align*}
Combining this with the following result gives Theorem \ref{maintheorem-for-characters}.
\begin{theorem} Assume that $D^{i}(z,s)$ is regular at $s_j=1/2+ir_j$ for $i=0,\ldots,n-1.$ Then
  \begin{enumerate}
\item the function $\varphi^{(l)}(s)$ is regular at $s_j$ for $l=0,1,\ldots, 2n-1$.
\item the function     $\varphi^{(2n)}(s)$ has at most a simple pole at $s_j$. Furthermore the residue at $s_j$ is given by
\begin{equation*}
  \res_{s=s_j}\varphi^{(2n)}(s)=-\varphi(s_j)\binom{2n}{n}\norm{\res_{s=s_j}D^{n}(z,s)}^2.
\end{equation*}
  \end{enumerate}
\end{theorem}
This theorem is proved by investigating further \eqref{phi-n-expr} and \eqref{laplacian-to-D}. We refer to \cite[Thm 3.3]{PetridisRisager:2013a} for details.

\bibliographystyle{alpha}

\begin{bibdiv}
\begin{biblist}

\bib{Avelin:2007b}{article}{
      author={Avelin, Helen},
       title={Deformation of $\Gamma_0(5)$-cusp forms},
        date={2007},
        ISSN={0025-5718},
     journal={Math. Comp.},
      volume={76},
      number={257},
       pages={361\ndash 384},
      review={\MR{2261026}},
}

\bib{Balslev:1997a}{article}{
      author={Balslev, Erik},
       title={Spectral deformation of {{Laplacians}} on hyperbolic manifolds},
        date={1997},
        ISSN={1019-8385},
     journal={Comm. Anal. Geom.},
      volume={5},
      number={2},
       pages={213\ndash 247},
      review={\MR{MR1483980 (99j:58213)}},
}

\bib{BalslevCombes:1971a}{article}{
      author={Balslev, Erik},
      author={Combes, Jean-Michel},
       title={Spectral properties of many-body {{Schr{\"o}dinger}} operators
  with dilatation- analytic interactions.},
    language={English},
        date={1971},
     journal={Commun. Math. Phys.},
      volume={22},
       pages={280\ndash 294},
}

\bib{BruggemanDiamantis:2012a}{article}{
      author={Bruggeman, Roelof},
      author={Diamantis, Nikolaos},
       title={Higher-order {{Maass}} forms},
        date={2012},
        ISSN={1937-0652},
     journal={Algebra Number Theory},
      volume={6},
      number={7},
       pages={1409\ndash 1458},
      review={\MR{3007154}},
}

\bib{Buser:1992a}{book}{
      author={Buser, Peter},
       title={Geometry and spectra of compact {{Riemann}} surfaces},
      series={Progress in {{Mathematics}}},
   publisher={{Birkh{\"a}user Boston Inc.}},
     address={{Boston, MA}},
        date={1992},
      volume={106},
        ISBN={0-8176-3406-1},
      review={\MR{MR1183224 (93g:58149)}},
}

\bib{BalslevVenkov:2000}{incollection}{
      author={Balslev, Erik},
      author={Venkov, Alexei},
       title={Selberg's eigenvalue conjecture and the {{Siegel}} zeros for
  {{Hecke}} ${{L}}$-series},
        date={2000},
   booktitle={Analysis on homogeneous spaces and representation theory of
  {{Lie}} groups, {{Okayama}}\textendash{{Kyoto}} (1997)},
      series={Adv. {{Stud}}. {{Pure Math}}.},
      volume={26},
   publisher={{Math. Soc. Japan, Tokyo}},
       pages={19\ndash 32},
}

\bib{BalslevVenkov:2001a}{article}{
      author={Balslev, Erik},
      author={Venkov, Alexei},
       title={Spectral theory of {{Laplacians}} for {{Hecke}} groups with
  primitive character},
        date={2001},
        ISSN={0001-5962},
     journal={Acta Math.},
      volume={186},
      number={2},
       pages={155\ndash 217},
      review={\MR{MR1846029 (2002f:11057)}},
}

\bib{BalslevVenkov:2005a}{article}{
      author={Balslev, Erik},
      author={Venkov, Alexei},
       title={On the relative distribution of eigenvalues of exceptional
  {{Hecke}} operators and automorphic {{Laplacians}}},
        date={2005},
        ISSN={0234-0852},
     journal={Algebra i Analiz},
      volume={17},
      number={1},
       pages={5\ndash 52},
      review={\MR{MR2140673 (2007b:11071)}},
}

\bib{BalslevVenkov:2007a}{inproceedings}{
      author={Balslev, Erik},
      author={Venkov, A.~B.},
       title={Perturbation of embedded eigenvalues of {{Laplacians}}},
        date={2007},
   booktitle={Traces in {{Number Theory}}, {{Geometry}} and {{Quantum
  Fields}}},
      editor={Albeverio, Sergio},
      editor={Marcolli, Matilde},
      editor={Paycha, Sylvie},
      editor={Plazas, Jorge},
      volume={38},
   publisher={{Vieweg Verlag}},
       pages={ix, 223 pp},
}

\bib{BalslevVenkov:1998a}{article}{
      author={Balslev, Erik},
      author={Venkov, Alexei~B.},
       title={The {{Weyl}} law for subgroups of the modular group},
        date={1998},
        ISSN={1016-443X},
     journal={Geom. Funct. Anal.},
      volume={8},
      number={3},
       pages={437\ndash 465},
      review={\MR{MR1631251 (99g:11069)}},
}

\bib{DeshouillersIwaniecPhillipsSarnak:1985a}{article}{
      author={Deshouillers, Jean-Marc},
      author={Iwaniec, Henryk},
      author={Phillips, Ralph~S.},
      author={Sarnak, Peter},
       title={Maass {{Cusp Forms}}},
        date={1985},
     journal={Proceedings of the National Academy of Sciences},
      volume={82},
      number={11},
       pages={3533\ndash 3534},
}

\bib{FarmerLemurell:2005a}{article}{
      author={Farmer, David~W.},
      author={Lemurell, Stefan},
       title={Deformations of {{Maass}} forms},
        date={2005},
     journal={Mathematics of Computation},
      volume={74},
       pages={1967\ndash 1982},
}

\bib{Goldfeld:1999b}{incollection}{
      author={Goldfeld, Dorian},
       title={The distribution of modular symbols},
        date={1999},
   booktitle={Number theory in progress, {{Vol}}. 2
  ({{Zakopane}}-{{Ko}}\textbackslash{}'scielisko, 1997)},
   publisher={{de Gruyter}},
     address={{Berlin}},
       pages={849\ndash 865},
      review={\MR{2000g:11040}},
}

\bib{Goldfeld:1999a}{incollection}{
      author={Goldfeld, Dorian},
       title={Zeta functions formed with modular symbols},
        date={1999},
   booktitle={Automorphic forms, automorphic representations, and arithmetic
  ({{Fort Worth}}, {{TX}}, 1996)},
      series={Proc. {{Sympos}}. {{Pure Math}}.},
      volume={66},
   publisher={{Amer. Math. Soc.}},
     address={{Providence, RI}},
       pages={111\ndash 121},
      review={\MR{2000g:11039}},
}

\bib{Hecke:1936}{article}{
      author={Hecke, Erich},
       title={{\"U}ber die {{Bestimmung Dirichletscher Reihen}} durch ihre
  {{Funktionalgleichung}}},
        date={1936},
        ISSN={0025-5831},
     journal={Math. Ann.},
      volume={112},
      number={1},
       pages={664\ndash 699},
      review={\MR{MR1513069}},
}

\bib{Hejhal:1983a}{book}{
      author={Hejhal, Dennis~A.},
       title={The {{Selberg}} trace formula for {{PSL}}(2,\, {{R}}). {{Vol}}.
  2},
      series={Lecture {{Notes}} in {{Mathematics}}},
   publisher={{Springer-Verlag}},
     address={{Berlin}},
        date={1983},
      volume={1001},
        ISBN={3-540-12323-7},
      review={\MR{MR711197 (86e:11040)}},
}

\bib{Ivrii:2016c}{article}{
      author={Ivrii, Victor},
       title={100 years of {{Weyl}}'s law},
        date={2016},
        ISSN={1664-3607},
     journal={Bull. Math. Sci.},
      volume={6},
      number={3},
       pages={379\ndash 452},
      review={\MR{3556544}},
}

\bib{Iwaniec:2002a}{book}{
      author={Iwaniec, Henryk},
       title={Spectral methods of automorphic forms},
     edition={Second},
      series={Graduate {{Studies}} in {{Mathematics}}},
   publisher={{American Mathematical Society}},
     address={{Providence, RI}},
        date={2002},
      volume={53},
        ISBN={0-8218-3160-7},
      review={\MR{MR1942691 (2003k:11085)}},
}

\bib{JorgensonOSullivan:2008a}{article}{
      author={Jorgenson, Jay},
      author={O'Sullivan, Cormac},
       title={Unipotent vector bundles and higher-order non-holomorphic
  {{Eisenstein}} series},
        date={2008},
        ISSN={1246-7405},
     journal={J. Th{\'e}or. Nombres Bordeaux},
      volume={20},
      number={1},
       pages={131\ndash 163},
      review={\MR{MR2434161}},
}

\bib{Kato:1976a}{book}{
      author={Kato, Tosio},
       title={Perturbation theory for linear operators},
     edition={Second},
   publisher={{Springer-Verlag}},
     address={{Berlin}},
        date={1976},
        note={Grundlehren der Mathematischen Wissenschaften, Band 132},
      review={\MR{MR0407617 (53,11389)}},
}

\bib{LaxPhillips:1976a}{book}{
      author={Lax, Peter~D.},
      author={Phillips, Ralph~S.},
       title={Scattering theory for automorphic functions},
   publisher={{Princeton Univ. Press}},
     address={{Princeton, N.J.}},
        date={1976},
        ISBN={0-691-08179-4},
        note={Annals of Mathematics Studies, No. 87},
      review={\MR{MR0562288 (58,27768)}},
}

\bib{Luo:2001a}{article}{
      author={Luo, Wenzhi},
       title={Nonvanishing of ${{L}}$-values and the {{Weyl}} law},
        date={2001},
        ISSN={0003-486X},
     journal={Ann. of Math. (2)},
      volume={154},
      number={2},
       pages={477\ndash 502},
      review={\MR{MR1865978 (2002i:11084)}},
}

\bib{Maass:1949a}{article}{
      author={Maass, Hans},
       title={{\"U}ber eine neue {{Art}} von nichtanalytischen automorphen
  {{Funktionen}} und die {{Bestimmung Dirichletscher Reihen}} durch
  {{Funktionalgleichungen}}},
        date={1949},
        ISSN={0025-5831},
     journal={Math. Ann.},
      volume={121},
       pages={141\ndash 183},
      review={\MR{MR0031519 (11,163c)}},
}

\bib{MinakshisundaramPleijel:1949a}{article}{
      author={Minakshisundaram, Subbaramiah},
      author={Pleijel, Åke},
       title={Some {{properties}} of the {{eigenfunctions}} of {{the
  Laplace}}-{{operator}} on {{Riemannian manifolds}}},
    language={en},
        date={1949},
        ISSN={0008-414X, 1496-4279},
     journal={Canadian Journal of Mathematics},
      volume={1},
      number={3},
       pages={242\ndash 256},
}

\bib{Muller:1992a}{article}{
      author={M{\"u}ller, Werner},
       title={Spectral geometry and scattering theory for certain complete
  surfaces of finite volume},
        date={1992},
        ISSN={0020-9910},
     journal={Invent. Math.},
      volume={109},
      number={2},
       pages={265\ndash 305},
      review={\MR{93g:58151}},
}

\bib{Petridis:2002a}{article}{
      author={Petridis, Yiannis~N.},
       title={Spectral deformations and {{Eisenstein}} series associated with
  modular symbols},
        date={2002},
        ISSN={1073-7928},
     journal={Int. Math. Res. Not.},
      number={19},
       pages={991\ndash 1006},
      review={\MR{1 903 327}},
}

\bib{Petridis:1994c}{article}{
      author={Petridis, Yiannis~N.},
       title={On the singular set, the resolvent and {{Fermi}}'s golden rule
  for finite volume hyperbolic surfaces},
        date={1994},
        ISSN={0025-2611},
     journal={Manuscripta Math.},
      volume={82},
      number={3-4},
       pages={331\ndash 347},
      review={\MR{MR1265004 (95b:11053)}},
}

\bib{Petridis:1994b}{article}{
      author={Petridis, Yiannis~N.},
       title={Spectral data for finite volume hyperbolic surfaces at the bottom
  of the continuous spectrum},
        date={1994},
        ISSN={0022-1236},
     journal={J. Funct. Anal.},
      volume={124},
      number={1},
       pages={61\ndash 94},
      review={\MR{MR1284603 (95g:11045)}},
}

\bib{PetridisRisager:2004a}{article}{
      author={Petridis, Yiannis~N.},
      author={Risager, Morten~S.},
       title={Modular symbols have a normal distribution},
        date={2004},
        ISSN={1016-443X},
     journal={Geom. Funct. Anal.},
      volume={14},
      number={5},
       pages={1013\ndash 1043},
      review={\MR{MR2105951 (2005h:11101)}},
}

\bib{PetridisRisager:2013a}{article}{
      author={Petridis, Yiannis~N.},
      author={Risager, Morten~S.},
       title={Dissolving of cusp forms: {{Higher}} order {{Fermi}}'s golden
  rule},
        date={2013},
     journal={Mathematika},
      volume={59},
      number={2},
       pages={269\ndash 301},
}

\bib{PhillipsSarnak:1985a}{article}{
      author={Phillips, Ralph~S.},
      author={Sarnak, Peter},
       title={On cusp forms for co-finite subgroups of {{PSL}}(2, {{R}})},
        date={1985},
        ISSN={0020-9910},
     journal={Invent. Math.},
      volume={80},
      number={2},
       pages={339\ndash 364},
      review={\MR{MR788414 (86m:11037)}},
}

\bib{PhillipsSarnak:1985c}{article}{
      author={Phillips, Ralph~S.},
      author={Sarnak, Peter},
       title={The {{Weyl}} theorem and the deformation of discrete groups},
        date={1985},
        ISSN={0010-3640},
     journal={Comm. Pure Appl. Math.},
      volume={38},
      number={6},
       pages={853\ndash 866},
      review={\MR{MR812352 (87f:11035)}},
}

\bib{PhillipsSarnak:1992a}{article}{
      author={Phillips, Ralph~S.},
      author={Sarnak, Peter},
       title={Perturbation theory for the {{Laplacian}} on automorphic
  functions},
        date={1992},
        ISSN={0894-0347},
     journal={J. Amer. Math. Soc.},
      volume={5},
      number={1},
       pages={1\ndash 32},
      review={\MR{MR1127079 (92g:11056)}},
}

\bib{Roelcke:1953a}{article}{
      author={Roelcke, Walter},
       title={{\"U}ber die {{Wellengleichung}} bei {{Grenzkreisgruppen}} erster
  {{Art}}},
        date={1953},
     journal={S.-B. Heidelberger Akad. Wiss. Math.-Nat. Kl.},
      volume={1953/1955},
       pages={159\ndash 267 (1956)},
      review={\MR{18,476d}},
}

\bib{Roelcke:1966a}{article}{
      author={Roelcke, Walter},
       title={Das {{Eigenwertproblem}} der automorphen {{Formen}} in der
  hyperbolischen {{Ebene}}. {{I}}, {{II}}.},
         how={Math. Ann. 167 (1966), 292--337; ibid 168 (1966), 261--324.},
        date={1966},
}

\bib{Sarnak:1990b}{incollection}{
      author={Sarnak, Peter},
       title={On cusp forms. {{II}}},
        date={1990},
   booktitle={Festschrift in honor of {{I}}. {{I}}. {{Piatetski}}-{{Shapiro}}
  on the occasion of his sixtieth birthday, {{Part II}} ({{Ramat Aviv}},
  1989)},
      series={Israel {{Math}}. {{Conf}}. {{Proc}}.},
      volume={3},
   publisher={{Weizmann}},
     address={{Jerusalem}},
       pages={237\ndash 250},
      review={\MR{MR1159118 (93e:11068)}},
}

\bib{Selberg:1956a}{article}{
      author={Selberg, Atle},
       title={Harmonic analysis and discontinuous groups in weakly symmetric
  {{Riemannian}} spaces with applications to {{Dirichlet}} series},
        date={1956},
     journal={J. Indian Math. Soc. (N.S.)},
      volume={20},
       pages={47\ndash 87},
      review={\MR{MR0088511 (19,531g)}},
}

\bib{Selberg:1989a}{book}{
      author={Selberg, Atle},
       title={G{\"o}ttingen lecture notes in {{Collected}} papers. {{Vol}}.
  {{I}}},
   publisher={{Springer-Verlag}},
     address={{Berlin}},
        date={1989},
        ISBN={3-540-18389-2},
        note={With a foreword by K. Chandrasekharan},
      review={\MR{MR1117906 (92h:01083)}},
}

\bib{Simon:1973a}{article}{
      author={Simon, Barry},
       title={Resonances in n-body quantum systems with dilatation analytic
  potentials and the foundations of time-dependent perturbation theory},
        date={1973},
        ISSN={0003-486X},
     journal={Ann. of Math. (2)},
      volume={97},
       pages={247\ndash 274},
      review={\MR{MR0353896 (50,6378)}},
}

\bib{Venkov:1979}{article}{
      author={Venkov, Alexei~B.},
       title={Artin-{{Takagi}} formula for the {{Selberg}} zeta function and
  the {{Roelcke}} hypothesis},
        date={1979},
        ISSN={0002-3264},
     journal={Doklady Akademii Nauk SSSR},
      volume={247},
      number={3},
       pages={540\ndash 543},
}

\bib{Venkov:1990a}{article}{
      author={Venkov, Alexei~B.},
       title={A remark on the discrete spectrum of the automorphic laplacian
  for a generalized cycloidal subgroup of a general fuchsian group},
        date={1990},
     journal={Journal of Mathematical Sciences},
      volume={52},
      number={3},
       pages={3016\ndash 3021},
}

\bib{Venkov:1990c}{book}{
      author={Venkov, Alexei~B.},
       title={Spectral theory of automorphic functions and its applications},
      series={Mathematics and Its {{Applications}} ({{Soviet Series}})},
   publisher={{Kluwer Academic Publishers Group}},
     address={{Dordrecht}},
        date={1990},
      volume={51},
        ISBN={0-7923-0487-X},
        note={Translated from the Russian by N. B. Lebedinskaya},
      review={\MR{MR1135112 (93a:11046)}},
}

\bib{Venkov:1992}{article}{
      author={Venkov, Alexei~B.},
       title={On a multidimensional variant of the {{Roelcke}}-{{Selberg}}
  conjecture},
        date={1992},
        ISSN={0234-0852},
     journal={Rossi{\u \i}skaya Akademiya Nauk. Algebra i Analiz},
      volume={4},
      number={3},
       pages={145\ndash 158},
}

\bib{Weyl:1911}{article}{
      author={Weyl, Hermann},
       title={{Ueber die asymptotische Verteilung der Eigenwerte}},
        date={1911},
     journal={Nachrichten von der Gesellschaft der Wissenschaften zu
  G{\"o}ttingen, Mathematisch-Physikalische Klasse},
      volume={1911},
       pages={110\ndash 117},
}

\end{biblist}
\end{bibdiv}

\end{document}